\documentclass[12pt]{article}

\usepackage{latexsym,amsmath}
\usepackage{amssymb}
\usepackage{microtype}

\newtheorem{thm}{Theorem}[section]

\newtheorem{lemma}[thm]{Lemma}

\newtheorem{cor}[thm]{Corollary}

\newcommand{\beq}[1]{\begin{equation}\label{#1}}
\newcommand{\enq}[0]{\end{equation}}

\newcommand{\qed}[0]{{\hspace*{\fill}\mbox{$\Box$}}}

\newcommand{\cA}[0]{{\cal A}}

\newcommand{\cE}[0]{{\cal E}}

\newcommand{\cG}[0]{{\cal G}}

\newcommand{\cI}[0]{{\cal I}}
\newcommand{\cJ}[0]{{\cal J}}

\newcommand{\cO}[0]{{\cal O}}

\newcommand{\cW}[0]{{\cal W}}

\newcommand{\Z}{{\mathbb Z}}

\newcommand{\gd}[0]{\delta}

\newcommand{\grg}[0]{\gamma}

\newcommand{\gl}[0]{\lambda}

\newcommand{\gS}[0]{\Sigma}

\newcommand{\bestgl}[0]{\frac{c\log d}{d^{1/3}}}
\newcommand{\gld}[0]{\frac{\gl}{2}\left(\frac{2}{1+\gl}\right)^d}

\begin{document}
\renewcommand{\thefootnote}{\fnsymbol{footnote}}

\title{A threshold phenomenon for random independent sets in the discrete hypercube}

\author{David Galvin\thanks{Department of Mathematics,
University of Notre Dame, 255 Hurley Hall, Notre Dame IN
46556; dgalvin1@nd.edu. Research supported in part by National Security Agency grant H98230-10-1-0364.}}

\date{\today}

\maketitle

\begin{abstract}
Let $I$ be an independent set drawn from the discrete $d$-dimensional hypercube $Q_d=\{0,1\}^d$ according to the hard-core distribution with parameter $\gl>0$ (that is, the distribution in which each independent set $I$ is chosen with probability proportional to $\gl^{|I|}$).  We show a sharp transition around $\gl=1$ in the appearance of $I$: for $\gl>1$, $\min\{|I \cap \cE|, |I \cap \cO|\}=0$ asymptotically almost surely, where $\cE$ and $\cO$ are the bipartition classes of $Q_d$, whereas for $\gl<1$, $\min\{|I \cap \cE|, |I \cap \cO|\}$ is asymptotically almost surely exponential in $d$.
The transition occurs in an interval whose length is of order $1/d$.

A key step in the proof is an estimation of $Z_\gl(Q_d)$, the sum over independent sets in $Q_d$ with each set $I$ given weight $\gl^{|I|}$ (a.k.a. the hard-core partition function). We obtain the asymptotics of $Z_\gl(Q_d)$ for $\gl>\sqrt{2}-1$, and nearly matching upper and lower bounds for $\gl \leq \sqrt{2}-1$, extending work of Korshunov and Sapozhenko. These bounds allow us to read off some very specific information about the structure of an independent set drawn according to the hard-core distribution.

We also derive a long-range influence result. For all fixed $\gl>0$, if $I$ is chosen from the independent sets of $Q_d$ according to the hard-core distribution with parameter $\gl$, conditioned on a particular $v \in \cE$ being in $I$, then the probability that another vertex $w$ is in $I$ is $o(1)$ for $w \in \cO$ but $\Omega(1)$ for $w \in \cE$.
\end{abstract}

\section{Introduction and statement of results}

The focus of this paper is the discrete hypercube $Q_d$. This is the
graph on vertex set $V=\{0,1\}^d$ with two strings adjacent if they
differ on exactly one coordinate. It is a $d$-regular bipartite
graph with bipartition classes $\cE$ and $\cO$, where $\cE$ is the
set of vertices with an even number of $1$'s. Note that
$|\cE|=|\cO|=2^{d-1}$. (For graph theory basics, see e.g.
\cite{MGT}).

An {\em independent set} in $Q_d$ is a set of vertices no two of
which are adjacent. Write $\cI(Q_d)$ for the set of independent sets
in $Q_d$. The {\em hard-core model} with parameter $\gl$ on $Q_d$ (abbreviated ${\rm hc}(\gl)$) is the probability distribution on $\cI(Q_d)$ in which each $I$ is chosen with probability proportional to $\gl^{|I|}$. This fundamental statistical physics model interprets vertices of a graph (in this case, $Q_d$) as sites that may or may not be occupied by massive particles, and edges as bonds between pairs of sites (encoding, for example, spatial proximity). The occupation rule is that bonded sites may not be simultaneously occupied, so a legal configuration of particles corresponds to an independent set in the graph. In this context $\gl$ represents a density parameter, with larger $\gl$ favouring denser configurations.
(For an introduction to the hard-core model from a combinatorial perspective, see for example \cite{BrightwellWinkler}.)

In \cite{Kahn}, Kahn used entropy methods to make an extensive study of the hard-core model on $Q_d$ (and regular bipartite graphs in general) for fixed $\gl>0$. One of the main results is that an independent set from $Q_d$ chosen according ${\rm hc}(\gl)$ exhibits {\em phase coexistence} --- it comes either predominantly from $\cE$ or predominantly from $\cO$. Specifically, it is shown in \cite{Kahn} that for fixed $\gl, \varepsilon>0$, and for $I$ chosen from $\cI(Q_d)$ according to ${\rm hc}(\gl)$, both of
\begin{eqnarray*}
\left||I|-\frac{\gl}{1+\gl}2^{d-1}\right| & \leq & \frac{2^d}{d^{1-\varepsilon}}, \\
\min\{|I\cap \cE|, |I \cap \cO|\} & \leq & \frac{2^d}{d^{1/2-\varepsilon}}
\end{eqnarray*}
hold asymptotically almost surely ({\em a.a.s.}), that is, with probability tending to $1$ as $d \rightarrow \infty$. Informally, the work of \cite{Kahn} demonstrates that for all fixed $\gl>0$, ${\rm hc}(\gl)$ is close to $\frac{1}{2}\mu_{\cE}+\frac{1}{2}\mu_\cO$ where $\mu_\cE$ (or $\mu_\cO$) is a random subset of $\cE$ (or $\cO$) in which each vertex is chosen to be in the set independently with probability $\frac{\gl}{1+\gl}$. (This is just ${\rm hc}(\gl)$ on $\cE$ (or $\cO$).)

Kahn's estimates on $|I|$ and $\min\{|I\cap \cE|, |I \cap \cO|\}$ do not involve $\gl$. Here we are able to obtain more precise estimates that capture the dependence on $\gl$ and in particular show that $\gl=1$ is a critical value around which a transition occurs in the nature of the phase coexistence: for $\gl>1$, the smaller of $|I \cap \cE|$, $|I \cap \cO|$ is a.a.s. $0$, whereas for $\gl < 1$, it is a.a.s. exponential in $d$. Allowing $\gl$ to vary with $d$, we find that the transition between $\min\left\{|I \cap \cE|,|I \cap \cO|\right\}$ being a.a.s. $0$ and a.a.s. going to infinity with $d$ occurs in an interval of length order $1/d$.

To state our results precisely we consider four possible ranges of $\gl$:
\begin{eqnarray}
& \gl \geq 1 + \frac{\omega(1)}{d} & \label{lambda1} \\
& |\gl-1| \leq \frac{O(1)}{d} & \label{lambda2} \\
& \sqrt{2}-1 + \frac{(\sqrt{2}+\Omega(1))\log d}{d} \leq \gl \leq 1 - \frac{\omega(1)}{d} & \label{lambda3} \\
& \bestgl \leq \gl \leq \sqrt{2}-1 + \frac{(\sqrt{2}+o(1))\log d}{d} & \label{lambda4}
\end{eqnarray}
where $c>0$ is an absolute constant (that we do not explicitly compute).
Here and in what follows, $\omega(1)$ indicates a function of $d$ that tends to infinity as $d$ does; $o(1)$ a function that tends to $0$; $\Omega(1)$ a function that is eventually always greater than some constant greater than $0$; and $O(1)$ a function that is bounded above by a constant. All implied constants will be independent of $d$, all limiting statements are as $d \rightarrow \infty$, and where we are not taking a limit, we will always assume that $d$ is large enough to support our assertions. Unless otherwise indicated, all logarithms are to base $e$.
\begin{thm} \label{thm-thresh}
Let $I$ be chosen from $\cI(Q_d)$ according to ${\rm hc}(\gl)$.
\begin{enumerate}
\item
For $\gl$ satisfying (\ref{lambda1}), a.a.s.
\begin{equation} \label{eq-max123}
\left|\max\{|I\cap \cE|, |I \cap \cO|\} - \frac{\gl 2^{d-1}}{1+\gl}\right| \leq 2^{d/2}\sqrt{\log d}
\end{equation}
and
\begin{equation} \label{eq-thresh_gl1}
\min\{|I\cap \cE|, |I \cap \cO|\} = 0.
\end{equation}
\item
For $\gl$ satisfying (\ref{lambda2}), a.a.s. (\ref{eq-max123}) holds. If $\gl = 1+\frac{k+o(1)}{d}$ for some constant $k$ then for each $c \in {\mathbb N}$
\begin{equation} \label{eq-scaling}
\Pr\left(\min\{|I\cap \cE|, |I \cap \cO|\} = c \right) \sim \frac{\left(\frac{1}{2}e^{-k/2}\right)^c}{c!}\exp\left\{-\frac{1}{2}e^{-k/2}\right\}.
\end{equation}
\item
For $\gl$ satisfying (\ref{lambda3}), a.a.s. (\ref{eq-max123}) holds, as well as
\begin{equation} \label{eq-thresh_gl3}
\frac{\left|\min\{|I\cap \cE|, |I \cap \cO|\} - \gld\right|}{\sqrt{(2+\varepsilon)\gld \log \left(\gld\right)}} \leq 1
\end{equation}
where $\varepsilon > 0$ is arbitrary.
\item
For $\gl$ satisfying (\ref{lambda4}), a.a.s.
\begin{equation} \label{eq-max4}
\left|\max\{|I\cap \cE|, |I \cap \cO|\} - \frac{\gl 2^{d-1}}{1+\gl}\right| \leq d(\log d)\left(\frac{2}{1+\gl}\right)^d
\end{equation}
and
\begin{equation} \label{eq-thresh_gl4_detail}
\frac{1}{4\log m}\gld \leq \min\{|I\cap \cE|, |I \cap \cO|\} \leq e m^2 \gld
\end{equation}
where $m=m(\gl,d) < d/\sqrt{\log d}$ is any natural number satisfying
\begin{equation} \label{def-mgl}
(ed^2)^m\gl^{m+1}(1+\gl)^{2m(m+1)}\frac{2^d}{(1+\gl)^{d(m+1)}} = o(1).
\end{equation}
\end{enumerate}
\end{thm}
The upper bound on $m$ in (\ref{def-mgl}) helps make our analysis more tractable, and does not impose any serious restriction: for any $\gl = \omega(\sqrt{\log d}/d)$, for example, $m$ can be taken to be $o(d/\sqrt{\log d})$.

\medskip

The following corollary of Theorem \ref{thm-thresh} is immediate.
\begin{cor}
$\Pr\left(\min\{|I\cap \cE|,|I \cap \cO|\}=0\right)$ goes from $1-o(1)$ to $o(1)$ as $\gl$ goes from $1+\omega(1/d)$ to $1-\omega(1/d)$.
\end{cor}
For fixed $\gl \leq \sqrt{2}-1$ we satisfy (\ref{def-mgl}) by taking $m = [1/\log_2(1+\gl)]$ and so combining (\ref{eq-thresh_gl3}) and (\ref{eq-thresh_gl4_detail}) we also get the following corollary.
\begin{cor}
For all fixed $\gl < 1$, there are constants $c(\gl)$ and $C(\gl)$ (independent of $d$) such that a.a.s.
$$
c(\gl) \left(\frac{2}{1+\gl}\right)^d \leq \min\{|I\cap \cE|, |I \cap \cO|\} \leq C(\gl) \left(\frac{2}{1+\gl}\right)^d.
$$
\end{cor}

Our proof of Theorem \ref{thm-thresh} provides structural information about the smaller of $I\cap \cE$ and $I \cap \cO$ for $I$ chosen according to ${\rm hc}(\gl)$. For simplicity, we state the following result for fixed $\gl$ (not varying with $d$). Say that $A \subseteq \cE$ (or $\cO$) is {\em $2$-linked} if $A \cup N(A)$ induces a connected subgraph, and note that each $A \subseteq \cE$ (or $\cO$) can be partitioned into {\em $2$-components} --- maximal $2$-linked subsets. Write $I_{\rm min}$ for the smaller of $I \cap \cE$ and $I \cap \cO$ (with an arbitrary choice being made if $|I \cap \cE|=|I \cap \cO|$).
\begin{thm} \label{thm-structure}
Fix $\gl > 0$ and let $I$ be chosen from $\cI(Q_d)$ according to ${\rm hc}(\gl)$. The following statements are all true a.a.s..
\begin{enumerate}
\item If $\gl > 1$ then $I_{\rm min}=\emptyset$.
\item If $\gl = 1$ then $I_{\rm min}$ consists of $k$ $2$-components, each of size $1$, with $k$ being drawn from a Poisson distribution with parameter $1/2$.
\item If $1 > \gl > \sqrt{2}-1$, then $I_{\rm min}$ consists of $k$ $2$-components, each of size $1$, where $k$ satisfies
$$
\left|k-\gld\right| \leq \sqrt{(2+\varepsilon)\gld \log\left(\gld\right)}.
$$
for any $\varepsilon > 0$.
\item If $2^{1/m}-1 \geq \gl > 2^{1/(m+1)}-1$ for some integer $m \geq 2$, then $I_{\rm min}$ consists of $k$ $2$-components, each of size at most $m$,
where $k$ satisfies
$$
\frac{1}{4\log m}\frac{\gl}{2} \left(\frac{2}{1+\gl}\right)^{d} \leq k \leq em\gld.
$$
\end{enumerate}
\end{thm}
As will be seen in our proof of Theorem \ref{thm-structure}, the fourth statement above is also true for $\gl$ satisfying (\ref{lambda4}) as long as $m$ is chosen to satisfy (\ref{def-mgl}).

\medskip

Theorem \ref{thm-structure} suggests that a sequence of threshold phenomena occur for independent sets chosen from $Q_d$ according to ${\rm hc}(\gl)$. The one we exhibit is that as $\gl$ passes across $1$, $I_{\rm min}$ goes from being empty to consisting of exponentially many singleton $2$-components (and nothing else), with these $2$-components arriving (in a Poisson manner) in a window of width $1/d$. It is tempting to conjecture that for each $m \geq 2$, as $\gl$ passes across $2^{1/m}-1$, $I_{\rm min}$ goes from having exponentially many $2$-components of size $m-1$ (and smaller), and no $2$-components of size $m$, to having exponentially many $2$-components of size $m$ (and smaller), with these new $2$-components arriving (in an appropriate Poisson manner) in a short threshold window. An interesting direction for future work would be to extend what we have done for $\gl$ satisfying (\ref{lambda1}), (\ref{lambda2}) and (\ref{lambda3}), and determine, for $\gl$ satisfying (\ref{lambda4}) (and smaller), the exact number of $2$-components of size $k$ in $I_{\rm min}$ for each relevant $k$, and the exact nature of the transition across $2^{1/m}-1$ for each $m \geq 2$.

\medskip

To understand probabilities associated with the hard-core model, it is useful to understand the normalizing constant (or partition function)
$$
Z_\gl(Q_d) = \sum_{I \in \cI(Q_d)} \gl^{|I|}.
$$
In the case $\gl=1$, this just counts the number of independent sets in $Q_d$. Motivated by the interpretation of independent sets as ``codes of distance $2$'' over a binary alphabet, Korshunov and Sapozhenko \cite{KorshunovSapozhenko} gave an asymptotic estimate in this case.
\begin{thm} \label{thm-sapkor}
$|\cI(Q_d)| \sim 2 \sqrt{e}2^{2^{d-1}}$ as $d \rightarrow \infty$.
\end{thm}
The following theorem, which extends Theorem \ref{thm-sapkor} to a wider range of $\gl$, is the main tool in our approach to Theorem \ref{thm-thresh}.
\begin{thm} \label{thm-Zdest}
$$
Z_\gl(Q_d) = \left\{
\begin{array}{ll}
(2+o(1))(1+\gl)^{2^{d-1}} & \mbox{if $\gl$ satisfies (\ref{lambda1})} \\
(2+o(1))(1+\gl)^{2^{d-1}}\exp\left\{\gld\right\} & \mbox{if $\gl$ satisfies (\ref{lambda2}) or (\ref{lambda3})} \\
(1+\gl)^{2^{d-1}}\exp\left\{\gld(1+o(1))\right\} & \mbox{if $\gl$ satisfies (\ref{lambda4})}.
\end{array}
\right.
$$
\end{thm}
We expect that the range of validity for the third of these estimates can be extended to $\gl > (1+\Omega(1))\log d/d$. Indeed, we already know that the lower bound is valid for this range of $\gl$, since (\ref{int7}) and the second clause of (\ref{int16}), which combine to give the lower bound, both hold for $\gl$ in this range. On the other hand, the upper bound would follow immediately from an extension of Lemma \ref{lem-mainsaplemma_weighted} to the range $\gl > (1+\Omega(1))\log d/d$.

The gap between the upper and lower bounds for $Z_\gl(Q_d)$ for $\gl$ satisfying (\ref{lambda4}) is the main obstacle to understanding more precisely the evolution of the ${\rm hc}(\gl)$ independent set, as discussed after the statement of Theorem \ref{thm-structure}.

\medskip

The phenomenon of phase coexistence is, unsurprisingly, accompanied by a long-range influence phenomenon. An independent set $I$ chosen from
$\cI(Q_d)$ according to ${\rm hc}(\gl)$ is drawn (in the range of $\gl$ that we are considering) either predominantly from $\cE$ or predominantly from $\cO$. If we are given the information that a particular vertex ($v \in \cE$, say) is in $I$, then that should make it very likely that $I$ is drawn mostly from $\cE$. So if we then ask what is the probability that another vertex ($w$, say) is in $I$, the answer should depend on the parity of $w$, being quite small if $w \in \cO$ and reasonably large if $w \in \cE$. This heuristic can be made rigorous.
\begin{thm} \label{thm-influence}
Let $\gl$ satisfy $\gl > \bestgl$. Let $I$ be chosen from $\cI(Q_d)$ according to ${\rm hc}(\gl)$.
If $u,v \in \cE$ and $w \in \cO$ are three vertices in $Q_d$
then
\begin{equation} \label{pt_odd}
\Pr\left(u \in I~|~w \in I\right)
\leq (1+\gl)^{-d(1-o(1))}
\end{equation}
and
\begin{equation} \label{pt_even}
\Pr\left(u \in I~|~v \in I\right)  \geq
\frac{\gl}{1+\gl}(1-o(1)).
\end{equation}
\end{thm}

\medskip

Estimates of $Z_\gl(Q_d)$ can also be used to obtain information on the number of independent sets of $Q_d$ of a given size; this topic will be explored in detail in a subsequent paper \cite{Galvin4}.

\medskip

Very few specific properties of $Q_d$ are used in the sequel. We just use the fact that it is a regular bipartite graph which satisfies certain isoperimetric bounds (specifically, those of Lemmas \ref{lem_easy-cube-iso} and \ref{lem-main-cube-iso}). Analogues of all of our main theorems could be obtained for any family of regular bipartite graphs with appropriate isometric properties, but in the absence of an appealing general formulation,
we confine ourselves here to considering $Q_d$.

\medskip

An overview of our approach is given in Section \ref{sec-overview}. The main technical lemma (Lemma \ref{lem-mainsaplemma_weighted}) is stated in Section \ref{sec-NT}, along with notation and other useful lemmas. The proofs of all the stated theorems appear in Section \ref{sec-proofs}, and a proof of Lemma \ref{lem-mainsaplemma_weighted} is presented in Section \ref{sec-sap_proof}.

\section{Overview} \label{sec-overview}

A trivial lower bound on $Z_\gl(Q_d)$ for all $\gl>0$ is $2(1+\gl)^{2^{d-1}}-1$: just consider the contribution from those sets which are drawn either entirely from $\cE$ or entirely from $\cO$. To improve this to the lower bounds appearing in Theorem \ref{thm-Zdest}, we consider not just
independent sets which are confined purely to either $\cE$ or $\cO$.
It is easy to see that there is a contribution of
$$
2^{d-1} \gl(1+\gl)^{2^{d-1}-d}=(1+\gl)^{2^{d-1}}\gld
$$
from those independent sets that have just one vertex from $\cO$, (and the same from those that have just one vertex from $\cE$) and more
generally a contribution of approximately
$$
2(1+\gl)^{2^{d-1}} \frac{1}{k!} \left(\gld\right)^k
$$
from those independent sets which consist of exactly $k$ non-nearby vertices on one side of the bipartition, for reasonably small $k$ (by ``non-nearby'' it is meant that there are
no common neighbours between pairs of the vertices). Indeed, there
are $2$ ways to chose the bipartition class that has $k$ vertices, and approximately ${2^{d-1} \choose k} \approx \frac{1}{k!}2^{(d-1)k}$ ways to choose the $k$
vertices. These vertices together have a neighbourhood of
size $kd$, so the sum of the weights of independent sets that extend the $k$
vertices is $\gl^k(1+\gl)^{2^{d-1}-kd}$ .

Summing over $k$ we get a lower bound on $Z_\gl(Q_d)$ of approximately
$$
2(1+\gl)^{2^{d-1}} \exp\left\{\gld\right\}.
$$
This lower bound could have also been achieved by summing only from $k$ a little below to a little above $\gld$ (where the mass of the Taylor series of $\exp\left\{\gld\right\}$ is concentrated) and, once $k$ vertices have been chosen from one side, only considering extensions to the other side which have close to $\frac{\gl 2^{d-1}}{1+\gl}$ vertices (where the mass of the binomial series $(1+\gl)^{2^{d-1}}$ is concentrated). This does not cause the count of extensions to drop much below $(1+\gl)^{2^{d-1}-dk}$ as long as $dk$ is much smaller than $2^{d-1}$, which it will be for $k \approx \gld$ and $\gl > \bestgl$. In this way we see that the lower bound on $Z_\gl(Q_d)$ can be achieved by only considering independent sets $I$ with $\min\{|I \cap \cE|,|I \cap \cO|\} \approx \gld$ and $\max\{|I \cap \cE|,|I \cap \cO|\} \approx \frac{\gl 2^{d-1}}{1+\gl}$. Thus an upper bound that matches the lower bound completes the proofs of both Theorems \ref{thm-Zdest} and \ref{thm-thresh} (as well as Theorem \ref{thm-structure}, as we shall see).

\medskip

To motivate the upper bound, consider what happens when we count the contribution from
independent sets that have exactly two nearby vertices from $\cO$
(that is, two vertices with a common neighbour). There are
approximately $d^22^{d-2}$ choices for this pair (as opposed to
approximately $2^{2d-2}$ choices for a pair of vertices without a
common neighbour), since once the first vertex has been chosen the
second must come from the approximately $d^2/2$ vertices at distance
two from the first. The sum of the weights of independent sets that extend each choice is $\gl^2(1+\gl)^{2^{d-1}-2d+2}$, roughly the same as
the sum of the weights of extensions in the case of the pair of vertices without
a common neighbour. The key point here is that any pair of
vertices from $\cO$ has at most two neighbours in common, so has neighbourhood size of approximately $2d$,
whether or not the vertices are nearby. Thus we get an
additional contribution of approximately
$$
(1+\gl)^{2^{d-1}}\frac{d^22^{d-2}}{(1+\gl)^{2d}}
$$
to the partition function from those sets with two nearby vertices
from $\cO$, negligible compared to the addition contribution to the partition function from those sets
with two non-nearby vertices from $\cO$.

The main work in upper
bounding $Z_\gl(Q_d)$ involves extending this to the
observation that the only non-negligible contribution to the partition function
comes from independent sets that on one side consist of a set of
vertices with non-overlapping neighbourhoods. This in turn amounts
to showing that there is a negligible contribution from those
independent sets which are $2$-linked on one side.
This entails proving a technical lemma (Lemma \ref{lem-mainsaplemma_weighted}) bounding
the sum of the weights of $2$-linked subsets of $\cE$ of a given size whose neighbourhood
in $\cO$ is of a given size. This lemma is a weighted generalization of an enumeration result originally introduced by Sapozhenko in \cite{Sapozhenko} and used in \cite{Sapozhenko2} to simplify the original proof of Theorem \ref{thm-sapkor}. A weaker form of Lemma \ref{lem-mainsaplemma_weighted} is proved in \cite{GalvinTetali2} where it used to estimate the weighted sum of independent sets in $Q_d$ satisfying $|I\cap \cE|=|I\cap \cO|$.

\section{Preliminaries} \label{sec-NT}

Let $\gS=(V,E)$ be a finite graph.
For $A \subseteq V$ write $N(A)$
for the set of vertices outside $A$ that are neighbours of a vertex
in $A$, and set
$$
[A]=\{v \in V:N(\{v\})\subseteq N(A)\};
$$
note that if $A$ is an independent set then $A \subseteq [A]$. For bipartite $\gS$ with bipartition $X \cup Y$ say that $A \subseteq X$ (or $Y$) is {\em small} if $|[A]| \leq |X|/2$ (or $|Y|/2$).

Say that $A$ is {\em $k$-linked} if for every $u,v\in A$ there is
a sequence $u=u_0, u_1, \ldots, u_l=v$ in $A$ such that for each $i = 0, \ldots, l-1$, the length of the shortest path connecting $u_i$ and $u_{i+1}$ is at most $k$, or, equivalently, if $A$ is connected in the graph obtained from $G$ by joining all pairs of vertices at distance at most $k$ (this extends our earlier notion of $2$-linkedness). The following easy lemma is from \cite{Sapozhenko}.
\begin{lemma} \label{nearbysets}
If $A$ is $k$-linked, and $T \subseteq V$ is such that for each $u \in A$ the length of the shortest path connecting $u$ to $T$ is at most $\ell$ and for each $v \in T$ the length of the shortest path connecting $v$ to $A$ is at most $\ell$, then $T$ is $(k+2\ell)$-linked.
\end{lemma}
Note that for bipartite $\gS$ and $A \subseteq X$ (or $Y$),
if $A$ is $2$-linked then so is $[A]$. Any $A$ can be decomposed into its maximal $2$-linked subsets; we refer to these as
the {\em $2$-components} of $A$.

\medskip

The following lemma bounds the number of connected subsets of a graph; see \cite[Lemma 2.1]{GalvinKahn}. (The bound given in \cite{GalvinKahn} is $(e\Delta)^n$, but the proof easily gives the claimed improvement.)
\begin{lemma} \label{lem-counting-conn-subgraphs}
Let $\gS$ have maximum degree $\Delta$. The number of $n$-vertex subsets of
$V$ which contain a fixed vertex and induce a connected subgraph is
at most $(e\Delta)^{n-1}$.
\end{lemma}
We will use the following easy corollary
which follows from the fact that a $k$-linked subset of a $d$-regular graph is
connected in a graph with all degrees at most $(k+1)d^k$.
\begin{cor} \label{Tree}
Let $\gS$ be $d$-regular with $d \geq 2$. The number of $k$-linked subsets of $V$ of size
$n$ which contain a fixed vertex is at most $\exp\left\{3nk\log d\right\}$.
\end{cor}

The next lemma is a special case of a fundamental result due to
Lov\'asz \cite{Lovasz} and Stein \cite{Stein}. For bipartite $\gS$ with bipartition $X \cup Y$, say that $Y' \subseteq Y$
{\em covers} $X$ if each $x \in X$ has a neighbour in $Y'$.
\begin{lemma} \label{lovaszstein}
If $\gS$ as above satisfies $|N(x)| \geq a$ for each $x \in X$
and $|N(y)| \leq b$ for each $y \in Y$ then there is some $Y' \subseteq Y$ that covers $X$ and satisfies
$$
|Y'| \leq (|Y|/a)(1 + \log b).
$$
\end{lemma}

The following is a special case of Hoeffding's Inequality \cite{Hoeffding}.
\begin{lemma} \label{lem-Hoeffding}
For all $\gl>0$, $\gd > 0$ and $m \in {\mathbb N}$,
$$
\sum_{j=\left\lfloor m \left(\frac{\gl}{1+\gl}-\gd\right) \right\rfloor}^{\left\lceil m\left(\frac{\gl}{1+\gl}+\gd\right) \right\rceil} \gl^j{m \choose j} \geq \left(1-2\exp\left\{-2\gd^2m\right\}\right)(1+\gl)^m.
$$
\end{lemma}

We will need to compare the exponential function $e^x$ to truncates $e_D(x) = \sum_{k=0}^D \frac{x^k}{k!}$ of its Taylor series; the following will be sufficient.
\begin{lemma} \label{lem-exp}
For any $y \leq x < z$ with $y, z \in {\mathbb N}$,
\begin{equation} \label{eq-exp_lower}
e_y(x) \leq \exp\left\{y\log\frac{ex}{y}+\log(y+1)\right\}
\end{equation}
and
\begin{equation} \label{eq-exp_upper}
e^{x}-e_z(x) \leq \exp\left\{z\log\frac{ex}{z}+\log\left(\frac{x}{z-x}\right)\right\}.
\end{equation}
\end{lemma}

\medskip

\noindent {\em Proof}:
We have
$$
e_y(x) = \sum_{k=0}^y \frac{x^k}{k!}  \leq (y+1)\frac{x^y}{y!} \leq \exp\left\{y\log\left(\frac{ex}{y}\right) + \log (y+1)\right\}
$$
and
$$
e^{x} - e_z(x) = \sum_{k=z+1}^\infty \frac{x^k}{k!} \leq \frac{x^z}{z!} \sum_{k=1}^\infty \left(\frac{x}{z}\right)^k \leq \exp\left\{z\log\left(\frac{ex}{z}\right) + \log \left(\frac{x}{z-x}\right)\right\}
$$
in both cases using $n!\geq (n/e)^n$.
\qed

\begin{cor} \label{cor-exp}
Let $\{x_d\}_{d=1}^\infty$ be such that $x_d \rightarrow \infty$. With $\varepsilon_1=\sqrt{c_1\log x_d/x_d}$ and $\varepsilon_2=\sqrt{c_2\log x_d/x_d}$ where $c_1>2$ and $c_2>1$ are constants, we have
$$
e_{[(1+\varepsilon_2)x_d]}(x_d) - e_{[(1-\varepsilon_1)x_d]}(x_d)\sim e^{x_d}.
$$
\end{cor}

\medskip

\noindent {\em Proof}:
Note that the function $f(t)=(1+t)\log\left(\frac{e}{1+t}\right)$ has a local maximum at $t=0$ and for $t=o(1)$ satisfies $f(t)=1-\frac{t^2}{2}+o(t^2)$.
From (\ref{eq-exp_lower}) we have $e_{[(1-\varepsilon_1)x_d]}(x_d) = o(e^{x_d})$
and from (\ref{eq-exp_upper}) we have $e^{x_d} - e_{[(1+\varepsilon_2)x_d]}(x_d) =  o(e^{x_d})$.
\qed

\medskip

For the remainder of this section, we specialize to $\gS=Q_d$. We will need the following isoperimetric bounds for $A \subseteq \cE$ (or $\cO$) (see \cite[Lemma 6.2]{Galvin} for the first and \cite[Lemma 1.3]{KorshunovSapozhenko} for the second).

\begin{lemma} \label{lem_easy-cube-iso}
There is a constant $C_{\rm iso}>0$ such that for $A \subseteq \cE$ (or $\cO$), if $|A| \leq d^4$ then $|A| \leq C_{\rm iso}|N(A)|/d$. If $|A| \leq d/10$, then $|N(A)| \geq d|A|-2|A|(|A|-1)$.
\end{lemma}

\begin{lemma} \label{lem-main-cube-iso}
For $A \subseteq \cE$ (or $\cO$), if $|A| \leq 2^{d-2}$ then
$$
|N(A)| \geq \left(1+\Omega(1/\sqrt{d})\right)|A|.
$$
\end{lemma}

Our main tool is a weighted version of a result of Sapozhenko \cite{Sapozhenko}.
\begin{lemma} \label{lem-mainsaplemma_weighted}
For each $a, g \geq 1$, set
$$
{\cal G}(a,g) = \{A \subseteq \cE ~2{\rm -linked}:|[A]|=a
~\mbox{and}~|N(A)|=g\}.
$$
There are constants $c>0$ and $c'>0$ such that the following holds. If $\gl > \bestgl$ and $a \leq 2^{d-2}$
then
$$
\sum_{A \in {\cal G}(a,g)} \gl^{|A|} \leq
2^d(1+\gl)^g \exp\left\{-\frac{c'(g-a)\log d}{d^{2/3}}\right\}.
$$
\end{lemma}
Lemma \ref{lem-mainsaplemma_weighted} can be proved by combining \cite[Lemma 4.5]{Sapozhenko} with \cite[Lemmas 3.3 and 3.4]{GalvinTetali2} (indeed, a key result from \cite{GalvinTetali2} is a slightly weaker version of Lemma \ref{lem-mainsaplemma_weighted}). A proof of Lemma \ref{lem-mainsaplemma_weighted} is given in Section \ref{sec-sap_proof}. Here we establish the following corollary, which is all that we will use in the sequel.
\begin{cor} \label{cor-small_tail}
For $\gl > \bestgl$ and $m \leq d/\sqrt{\log d}$,
$$
\sum \gl^{|A|}(1+\gl)^{-|N(A)|} \leq (ed^2)^{m-1}\gl^m(1+\gl)^{2m(m-1)}\frac{2^d}{(1+\gl)^{md}}.
$$
where the sum is over all $A \subseteq \cE$ small and $2$-linked with $|A|\geq m$.
\end{cor}

\medskip

\noindent {\em Proof}: We consider the sum in three parts. Say that $A$ is {\em of type I} if
$|A| \leq d/10$; {\em of type II} if $d/10 < |A| \leq d^4$ and
{\em of type III} if $d^4 < |A|$.

For type I $A$ with $|A|=k$ ($k \geq m$) there are (by Lemma \ref{lem-counting-conn-subgraphs}) at most $2^{d-1}(ed^2)^{k-1}$ choices for $A$ (the factor of $2^{d-1}$ accounting for the choice of a fixed vertex in $A$ and the $d^2$ coming from the fact that each $A$ is connected in a graph with maximum degree at most $d^2$). By Lemma \ref{lem_easy-cube-iso} each such $A$ satisfies $|N(A)|\geq dk-2k(k-1)$. It follows that the contribution to the sum from type I $A$'s is at most
$$
\sum_{k = m}^{d/10} 2^{d-1}(ed^2)^{k-1} \gl^k(1+\gl)^{-dk+2k(k-1)}.
$$
For large enough $d$ (independent of $\gl$, in the range $\gl > \bestgl$) each summand above is at most one third its predecessor and so the total sum is at most
\begin{equation} \label{mainterm}
\frac{3}{4}(ed^2)^{m-1} \gl^{m}(1+\gl)^{2m(m-1)}\frac{2^d}{(1+\gl)^{md}}.
\end{equation}
To complete the proof of the corollary we will show that the contributions to the sum from $A$'s of type II and type III are negligible compared to (\ref{mainterm}). 

The contribution to the sum from type II $A$'s (again using Lemmas \ref{lem-counting-conn-subgraphs} and \ref{lem_easy-cube-iso}) is at most
$$
\sum_{k=d/10}^{d^4} 2^{d-1}(ed^2)^{k-1} \gl^k(1+\gl)^{-\frac{dk}{C_{\rm iso}}}
$$
(where $C_{\rm iso}$ is the constant from Lemma \ref{lem_easy-cube-iso}). For large enough $d$ (independent of $\gl$, in the range $\gl > \bestgl$) the first term in this sum is the largest and so the sum is at most
$$
d^4 2^{d-1}(ed^2)^{\frac{d}{10}-1} \gl^{\frac{d}{10}}(1+\gl)^{-\frac{d^2}{10C_{\rm iso}}}
$$
which is vanishingly small compared to (\ref{mainterm}) for all $\gl$ and $m$ in the specified range. 

In the range $|A|>d^4$ we partition the possible $A$'s according to $a:=|[A]|>d^4$ and $g:=|N(A)|>d^4$. By Lemma \ref{lem-mainsaplemma_weighted}, the sum over type III $A$'s is at most
$$
\sum_{a,g>d^4,~\cG(a,g)\neq \emptyset} \left(\sum_{A \in \cG(a,g)} \gl^{|A|}\right) (1+\gl)^{-g} \leq
\sum_{a,g>d^4} 2^d\exp\left\{-\frac{c'(g-a)\log d}{d^{2/3}}\right\}.
$$
By Lemma \ref{lem-main-cube-iso}, $g-a\geq \Omega(d^{7/2})$ and there are at most $2^{2d}$ choices for $a$ and $g$, so the sum is at most
$2^{-\omega(d^2)}$, which again is vanishingly small compared to (\ref{mainterm}) for all $\gl$ and $m$ in the specified range.
\qed

\section{Proofs of the main theorems}
\label{sec-proofs}

\subsection{Proof of Theorem \ref{thm-Zdest}}
\label{subsec-proof_of_thm_Zdest}

We will begin with a general upper bound on $Z_\gl(Q_d)$.
\begin{lemma} \label{lem-general_ub}
For any $\gl>0$,
$$
Z_\gl(Q_d) \leq 2 (1+\gl)^{2^{d-1}}
\exp\left\{\sum_{A \subseteq \cE~{\rm small},~2{\rm-linked},~|A|\geq 1}
\gl^{|A|}(1+\gl)^{-|N(A)|}\right\}.
$$
\end{lemma}

\medskip

To see that this implies the claimed upper bounds, note that
\begin{equation} \label{int10}
\sum  \gl^{|A|}(1+\gl)^{-|N(A)|} = \gld + \frac{\gl^2(1+\gl)^2}{4}{d \choose 2}\frac{2^d}{(1+\gl)^{2d}}
\end{equation}
where the sum is over all $A \subseteq \cE$ small and $2$-linked with $1\leq |A|\leq 2$.
The second term on the right corresponds to $|A|=2$: there are $2^{d-1}{d \choose 2}/2$ ways to choose $A \subseteq \cE$ small and $2$-linked with $|A|=2$, and each such $A$ has $|N(A)| = 2d -2$. The first term corresponds to $|A|=1$.
On the other hand, from Corollary \ref{cor-small_tail} we have that for all $\gl > \bestgl$
\begin{eqnarray}
\sum \gl^{|A|}(1+\gl)^{-|N(A)|} & \leq & (ed^2)^2\gl^3(1+\gl)^{12}\frac{2^d}{(1+\gl)^{3d}} \nonumber \\
& = & o\left(\frac{\gl^2(1+\gl)^2}{4}{d \choose 2}\frac{2^d}{(1+\gl)^{2d}}\right) \label{int11}
\end{eqnarray}
where the sum is now over all $A \subseteq \cE$ small and $2$-linked with $|A|\geq 3$.
Inserting (\ref{int10}) and (\ref{int11}) into Lemma \ref{lem-general_ub} we obtain (for $\gl > \bestgl$)
\begin{equation} \label{eq-part_func_ub}
Z_\gl(Q_d) \leq 2(1+\gl)^{2^{d-1}}\exp\left\{\gld + \gl^2(1+\gl)^2d^2\frac{2^d}{(1+\gl)^{2d}}\right\}.
\end{equation}

If $\gl=\gl(d)$ satisfies (\ref{lambda1})
then the exponent in (\ref{eq-part_func_ub}) is $o(1)$. If $\gl$ satisfies either (\ref{lambda2}) or (\ref{lambda3}) then
it is $\gld + o(1)$. Finally, if $\gl$ satisfies (\ref{lambda4}) then
it is $\gld(1 + o(1))$. This gives all the upper bounds of Theorem \ref{thm-Zdest}.

\medskip

\noindent {\em Proof of Lemma \ref{lem-general_ub}}:
A simple argument (based on the fact that $Q_d$
has a perfect matching) shows that for $I \in \cI(Q_d)$, at least one of $|[I \cap \cE]|\leq 2^{d-2}$, $|[I \cap \cO]|\leq 2^{d-2}$ holds.
By $\cE$-$\cO$ symmetry we therefore have
\begin{equation} \label{general_ub1}
Z_\gl(Q_d) \leq 2(1+\gl)^{2^{d-1}}\sum_{A \subseteq \cE~{\rm small}} \gl^{|A|}(1+\gl)^{-|N(A)|}.
\end{equation}
Decomposing $A$ into $2$-components $A_1, \ldots, A_k$, we have
\begin{equation} \label{eq-prod_form}
\gl^{|A|}(1+\gl)^{-|N(A)|} = \prod_{i=1}^k \gl^{|A_i|}(1+\gl)^{-|N(A_i)|}
\end{equation}
and
\begin{eqnarray}
\sum_{A \subseteq \cE~{\rm small}} \!\! \gl^{|A|}(1+\gl)^{-|N(A)|} & = &
\sum\left\{\prod_{i=1}^k \gl^{|A_i|}(1+\gl)^{-|N(A_i)|} : \!\!
\begin{array}{c}
k \geq 0 \\
A \subseteq \cE~{\rm small} \\
A=\cup_{i=1}^k A_i
\end{array}\right\}
 \nonumber \\
& \leq & \sum_{k\geq 0} \frac{\left(\sum
\gl^{|A|}(1+\gl)^{-|N(A)|}\right)^k}{k!}  \nonumber \\
& = & \exp\left\{\sum
\gl^{|A|}(1+\gl)^{-|N(A)|}\right\} \label{general_ub2}
\end{eqnarray}
where the unqualified sum in the last two lines is over all $A \subseteq \cE$ small and $2$-linked with $|A|\geq 1$.
Combining (\ref{general_ub2}) with (\ref{general_ub1}) we obtain the lemma. \qed

\medskip

Before turning to the lower bounds, we combine (\ref{int10}), (\ref{int11}) and (\ref{general_ub2}) to observe that
for $\gl > \frac{c\log d}{d^{1/3}}$ (for suitably large $c$)
we have
\begin{equation} \label{eq-inner_sum_upper}
\sum_{A \subseteq \cE~{\rm small}} \gl^{|A|}(1+\gl)^{-|N(A)|} \leq \exp\left\{\frac{\gl}{2}\left(\frac{2}{1+\gl}\right)^d + \frac{d^2 \gl^2 (1+\gl)^2 2^d}{(1+\gl)^{2d}}\right\}.
\end{equation}

\medskip

Now we turn to the lower bounds on $Z_\gl(Q_d)$, which will follow from a general bound that is more than what we need for the proof of Theorem \ref{thm-Zdest} but just what we need for much of Theorems \ref{thm-thresh} and \ref{thm-structure}.
\begin{lemma} \label{lem-general_lb}
For all $\gl \geq \frac{\omega(1)}{d}$ and $f \leq \ell \leq \frac{2^{d-2}}{d^2}$,
$$
Z_\gl(Q_d) \geq 2(1+\gl)^{2^{d-1}}\sum_{k=f}^{\ell}\frac{1}{k!}\left(\gld\right)^k
\exp\left\{-\frac{\ell^2 d^2}{2^{d-2}}\right\}\left(1-\frac{2}{d^2}\right).
$$
This lower bound is obtained by considering only $I$ which satisfy
$$
f \leq \min\{|I\cap \cE|, |I \cap \cO|\} \leq \ell
$$
and
$$
e_1 \leq \max\{|I\cap \cE|, |I \cap \cO|\} \leq e_2
$$
where
$$
e_1 = \frac{\gl}{1+\gl}\left(2^{d-1} - d\ell\right) - \sqrt{(\log d)\left(2^{d-1} - d f\right)}
$$
and
$$
e_2 = \frac{\gl}{1+\gl}\left(2^{d-1} - d f\right) + \sqrt{(\log d)\left(2^{d-1} - d f\right)}
$$
\end{lemma}

\medskip

Before proving the lemma, we use it to obtain the claimed lower bounds on $Z_\gl(Q_d)$ and complete the proof of Theorem \ref{thm-Zdest}.

For $\gl$ satisfying either (\ref{lambda1}) or (\ref{lambda2}), $\gld = O(1)$. With $f=f(\gl)=0$ and $\ell=\ell(\gl)=\log d$ (say), an application of (\ref{eq-exp_upper}) yields
$$
\sum_{k=f(\gl)}^{\ell(\gl)}\frac{1}{k!}\left(\gld\right)^k \geq \exp\left\{\gld\right\} -o(1)
$$
and we also have
$$
\exp\left\{-\frac{\ell(\gl)^2 d^2}{2^{d-2}}\right\} \geq 1-o(1).
$$
Putting these bounds into Lemma \ref{lem-general_lb} we get
\begin{equation} \label{int6}
Z_\gl(Q_d) \geq (2-o(1))(1+\gl)^{2^{d-1}} \exp\left\{\gld\right\}.
\end{equation}
Noting that $\gld = o(1)$ for $\gl$ satisfying (\ref{lambda1}), we get from (\ref{int6}) the claimed lower bounds on $Z_\gl(Q_d)$ for $\gl$ satisfying either (\ref{lambda1}) or (\ref{lambda2}).

For $\gl$ satisfying either (\ref{lambda3}) or (\ref{lambda4}), $\gld = \omega(1)$.
For any $\varepsilon>0$ set
$$
f(\gl) = \gld - \sqrt{(2+\varepsilon)\gld \log \left(\gld\right)}
$$
and
$$
\ell(\gl) = \gld + \sqrt{(2+\varepsilon)\gld \log \left(\gld\right)}.
$$
An application of Corollary \ref{cor-exp} yields
$$
\sum_{k=f(\gl)}^{\ell(\gl)}\frac{1}{k!}\left(\gld\right)^k \geq (1-o(1))\exp\left\{\gld\right\}
$$
and we also have
$$
\exp\left\{-\frac{\ell(\gl)^2 d^2}{2^{d-2}}\right\} \geq \exp\left\{-2d^2\gl^2\frac{2^d}{(1+\gl)^{2d}}\right\}.
$$
Putting these bounds into Lemma \ref{lem-general_lb} we get
\begin{equation} \label{int7}
Z_\gl(Q_d) \geq (2-o(1))(1+\gl)^{2^{d-1}} \exp\left\{\gld - 2d^2\gl^2\frac{2^d}{(1+\gl)^{2d}}\right\}.
\end{equation}
Noting that
\begin{equation} \label{int16}
2d^2\gl^2\frac{2^d}{(1+\gl)^{2d}} = \left\{
\begin{array}{ll}
o(1) & \mbox{for $\gl$ satisfying (\ref{lambda3})} \\
o\left(\gld\right) & \mbox{for $\gl$ satisfying (\ref{lambda4})}
\end{array}
\right.
\end{equation}
we get from (\ref{int7}) the claimed lower bounds on $Z_\gl(Q_d)$ for $\gl$ satisfying either (\ref{lambda3}) or (\ref{lambda4}).

\medskip

\noindent {\em Proof of Lemma \ref{lem-general_lb}}:
For each $f \leq k \leq \ell$, we consider the contribution to the partition function from those $I$ with $|I\cap \cE|=k$, $e_1(k) \leq |I \cap \cO| \leq e_2(k)$ and all $2$-components of $I \cap \cE$ having size $1$, where
$$
e_1(k):= \frac{\gl}{1+\gl}\left(2^{d-1} - dk\right) - \sqrt{(\log d)\left(2^{d-1} - dk\right)}
$$
and
$$
e_2(k):= \frac{\gl}{1+\gl}\left(2^{d-1} - dk\right) + \sqrt{(\log d)\left(2^{d-1} - dk\right)}.
$$
If we choose the elements of $I \cap \cE$ sequentially then each new vertex we add removes from consideration at most ${d \choose 2}+1\leq d^2$ vertices (those vertices which are at distance at most $2$ from the chosen vertex). So the number of choices for $I \cap \cE$ is at least
\begin{equation}  \label{counting_IintE}
\frac{\prod_{j=0}^{k-1}\left(2^{d-1}-jd^2\right)}{k!} \geq  \frac{2^{k(d-1)}}{k!}\left(1-\frac{\ell d^2}{2^{d-1}}\right)^{\ell}
\geq \frac{2^{k(d-1)}}{k!}\exp\left\{-\frac{\ell^2 d^2}{2^{d-2}}\right\},
\end{equation}
the second inequality using $1-x \geq e^{-2x}$ for $0 < x < 1/2$; the application is valid since $\ell \leq \frac{2^{d-2}}{d^2}$.

Once $I \cap \cE$ has been chosen, there are $2^{d-1}-dk$ vertices in $\cO$ from among which we choose between $e_1(k)$ and $e_2(k)$ to complete $I$. The sum of the weights of the valid extensions to $I$ is, using Lemma \ref{lem-Hoeffding},
\begin{equation} \label{counting_IintO}
\gl^{k} \sum_{j = e_1(k)}^{e_2(k)} \gl^j{2^{d-1}-dk \choose j}  \geq   \gl^{k}(1+\gl)^{2^{d-1}-dk}\left(1-\frac{2}{d^2}\right).
\end{equation}

Combining (\ref{counting_IintE}) and (\ref{counting_IintO}) and noting that $e_1 \leq e_1(k)$ and $e_2(k) \leq e_2$ for all $f \leq k \leq \ell$, we see that the contribution to the partition function from those $I$ with $f \leq |I\cap \cE| \leq \ell$ and $e_1 \leq |I \cap \cO| \leq e_2$ is at least
$$
(1+\gl)^{2^{d-1}}\sum_{k=f}^{\ell}\frac{1}{k!}\left(\gld\right)^k \exp\left\{-\frac{\ell^2 d^2}{2^{d-2}}\right\}\left(1-\frac{2}{d^2}\right).
$$
We get at least the same contribution from those $I$ with $f \leq |I\cap \cO| \leq \ell$, $e_1 \leq |I \cap \cE| \leq e_2$. Since $\ell < e_1$ there is no overlap between the two contributions, and all $I$ under consideration satisfy $f \leq \min\{|I\cap \cE|, |I \cap \cO|\} \leq \ell$ and $e_1 \leq \max\{|I\cap \cE|, |I \cap \cO|\} \leq e_2$. This completes the proof of the lemma.
\qed

\subsection{Proof of Theorem \ref{thm-thresh}}
\label{subsec-proof_of_thm_thresh}

The lower bounds on $Z_\gl(Q_d)$ for $\gl$ satisfying (\ref{lambda1}), (\ref{lambda2}) and (\ref{lambda3}) come from considering only $I$ satisfying
$$
b_1(\gl) \leq \max\left\{|I\cap \cE|, |I \cap\cO|\right\}-\frac{\gl 2^{d-1}}{1+\gl} \leq b_2(\gl)
$$
where
$$
b_1(\gl) = -d\ell(\gl) - \sqrt{(\log d)\left(2^{d-1}-df(\gl)\right)}
$$
and
$$
b_2 = - df(\gl) + \sqrt{(\log d)\left(2^{d-1}-df(\gl)\right)}
$$
with $f(\gl)$ and $\ell(\gl)$ as introduced in the discussion after the statement of Lemma \ref{lem-general_lb}. For all such $\gl$ we have
$$
b_1 \geq -2^{d/2}\sqrt{\log d}~~~\mbox{and}~~~b_2 \leq 2^{d/2}\sqrt{\log d},
$$
the main point in both cases being that for $\gl$ satisfying (\ref{lambda3}), $d\gl\left(\frac{2}{1+\gl}\right)^d=o(2^{d/2})$.
Since the lower bounds in this range are asymptotic to the upper bounds, that (\ref{eq-max123}) occurs a.a.s for this range of $\gl$ follows immediately, as does similarly the fact that (\ref{eq-thresh_gl3}) holds a.a.s. for $\gl$ satisfying (\ref{lambda3}).

That (\ref{eq-thresh_gl1}) holds a.a.s for $\gl$ satisfying (\ref{lambda1}) follows immediately from Theorem \ref{thm-Zdest}. Indeed, the contribution to $Z_\gl(Q_d)$ from those $I$ with $\min\{|I \cap \cE|, |I \cap \cO|\}=0$ is
$$
2(1+\gl)^{2^{d-1}}-1 \sim 2(1+\gl)^{2^{d-1}} \sim Z_\gl(Q_d).
$$

\medskip

We have to work a little harder to show that (\ref{eq-max4}) and (\ref{eq-thresh_gl4_detail}) occur a.a.s. for $\gl$ satisfying (\ref{lambda4}). In this range, set
$$
\cI_\cE(\gl) = \left\{I \in \cI(Q_d):\begin{array}{c}
{\rm cl}(I \cap \cE) \leq m\\
\frac{1}{4\log m}\frac{\gl}{2} \left(\frac{2}{1+\gl}\right)^{d} \leq k(I\cap \cE) \leq em\gld\\
\left|\max\{|I \cap \cE|,|I \cap \cO|\}-\frac{\gl 2^{d-1}}{1+\gl}\right| \leq d(\log d)\left(\frac{2}{1+\gl}\right)^d
\end{array}\right\}.
$$
(with $m$ as in (\ref{def-mgl})) where ${\rm cl}(A)$ and $k(A)$ are the size of the largest $2$-component of $A$ and the number of $2$-components of $A$, respectively, and define $\cI_\cO(\gl)$ analogously. Note that $\cI_\cE(\gl)$ and $\cI_\cO(\gl)$ are disjoint and that $I \in \cI_\cE(\gl)$ satisfies (\ref{eq-max4}) and (\ref{eq-thresh_gl4_detail}), so the following lemma completes the proof that (\ref{eq-max4}) and (\ref{eq-thresh_gl4_detail}) occur a.a.s. for $\gl$ satisfying (\ref{lambda4}).
\begin{lemma} \label{lem-small_no_of_comps_gl4}
For $\gl$ satisfying (\ref{lambda4}),
$$
Z_\gl(Q_d) \sim \sum_{I \in \cI_\cE(\gl)} \gl^{|I|} + \sum_{I \in \cI_\cO(\gl)} \gl^{|I|}.
$$
\end{lemma}

\medskip

\noindent {\em Proof of Lemma \ref{lem-small_no_of_comps_gl4}}: We begin by considering the contribution to $Z_\gl(Q_d)$ from those $I$ with $I \cap \cE$ small and ${\rm cl}(I \cap \cE)> m$. With the sum below over such $I$, and recalling (\ref{eq-prod_form}), we have
\begin{eqnarray}
\sum \gl^{|I|} & = & (1+\gl)^{2^{d-1}} \sum_{A \subseteq \cE ~{\rm small},~ {\rm cl}(A)> m} \gl^{|A|}(1+\gl)^{-|N(A)|} \nonumber \\
& \leq & (1+\gl)^{2^{d-1}}\sum_{A' \subseteq \cE~ {\rm small},~2-{\rm linked}, ~|A'|>m} \gl^{|A'|}(1+\gl)^{-|N(A')|} \times \nonumber \\
& & ~~~~~~~~~~~~~~~~~~~~~~~~~~~~~~~~~~~~~\sum_{A'' \subseteq \cE ~{\rm small}} \gl^{|A''|}(1+\gl)^{-|N(A'')|}
\nonumber \\
& \leq & Z_\gl(Q_d)\sum_{A' \subseteq \cE~ {\rm small},~2-{\rm linked}, ~|A'|> m} \gl^{|A'|}(1+\gl)^{-|N(A')|} \nonumber \\
& = & o\left(Z_\gl(Q_d)\right), \label{no-large-comps}
\end{eqnarray}
where in (\ref{no-large-comps}) we have used Corollary \ref{cor-small_tail}.
We similarly have a negligible contribution to $Z_\gl(Q_d)$ from those $I$ with $I \cap \cO$ small and ${\rm cl}(I \cap \cO)> m$.

Next we consider the contribution from those $I$ with $I \cap \cE$ small, ${\rm cl}(I \cap \cE) \leq m$ and $k(I \cap \cE) \leq \frac{1}{4\log m}\gld$. The contribution is at most
\begin{eqnarray}
& & (1+\gl)^{2^{d-1}}\sum_{k \leq \frac{1}{4\log m}\gld} {2^{d-1} \choose k} m^k \gl^k(1+\gl)^{-dk} \label{int14} \\
& \leq & (1+\gl)^{2^{d-1}} \sum_{k \leq \frac{1}{4\log m}\gld} \frac{m^k}{k!}\left(\gld\right)^k \nonumber \\
& \leq & (1+\gl)^{2^{d-1}} \exp\left\{(1-\Omega(1))\left(\gld\right)\right\} \label{int8} \\
& = & o\left(Z_\gl(Q_d)\right). \nonumber
\end{eqnarray}
The factor of ${2^{d-1} \choose k}$ in (\ref{int14}) counts the number of ways of choosing a fixed vertex in each of the $k$ $2$-components of $I \cap \cE$. The factor of $m^k$ counts the number of ways of assigning a size to each $2$-component. For each choice of a fixed vertex and a size ($\ell_i$, say) for each $2$-component, the contribution to $Z_\gl(Q_d)$ is at most
$$
(1+\gl)^{2^{d-1}}\prod_{i=1}^k(ed^2)^{\ell_i-1} \gl^{\ell_i}(1+\gl)^{-d\ell_i+2\ell_i(\ell_i-1)} \leq (1+\gl)^{2^{d-1}} \prod_{i=1}^k \gl (1+\gl)^{-d}
$$
(for large enough $d$, independent of $\gl$). In (\ref{int8}) we use (\ref{eq-exp_lower}).

A similar calculation (using (\ref{eq-exp_upper}) in place of (\ref{eq-exp_lower})) shows that the contribution from those $I$ with $I \cap \cE$ small, ${\rm cl}(I \cap \cE) \leq m$ and $k(I \cap \cE) \geq e m \gld$ is $o\left(Z_\gl(Q_d)\right)$, and by symmetry so too is the contribution from those $I$ with $I \cap \cO$ small, ${\rm cl}(I \cap \cO) \leq m$ and either $k(I \cap \cO) \leq \frac{1}{4\log m} \gld$ or $k \geq em \gld$.

We have shown that $\frac{1}{2}(1-o(1))$ of $Z_\gl(Q_d)$ comes from
$$
\cI'_\cE(\gl) = \left\{I \in \cI(Q_d): \begin{array}{c}
{\rm cl}(I \cap \cE) \leq m\\
\frac{1}{4\log m}\gld \leq k(I \cap \cE) \leq em \gld
\end{array}
\right\}
$$
and another $\frac{1}{2}(1-o(1))$ comes from the analogously defined $\cI'_\cO(\gl)$. (We have dropped  ``$I \cap \cE$ small'' since it is implied by the condition on $k(I \cap \cE)$).

It remains to show that the contribution to $\cI'_\cE(\gl)$ from those $I$ with $I \cap \cO$ either too large or too small is negligible. For each $\frac{1}{4\log m}\gld \leq k \leq em \gld$ and each choice of $k$ $2$-components $A_1, \ldots, A_k$ for $I \cap \cE$, the contribution to $\sum_{I \in \cI'_\cE(\gl)} \gl^{|I|}$ is
$$
\sum_{j=0}^{2^{d-1}-\sum_{i=1}^k |N(A_i)|} \gl^j {2^{d-1}-\sum_{i=1}^k |N(A_i)| \choose j}.
$$
By Lemma \ref{lem-Hoeffding}, all but a proportion at most $\frac{2}{d^2}$ of this sum comes from those $j$ satisfying
$$
\left| \frac{\gl2^{d-1}}{1+\gl}  - j \right| \leq  \frac{\gl}{1+\gl}\sum_{i=1}^k |N(A_i)|+\sqrt{(\log d) \left(2^{d-1}-\sum_{i=1}^k |N(A_i)|\right)}.
$$
For all $k$ in the range under consideration, and all possible choices of the $A_i$'s, we have
$$
\frac{\gl}{1+\gl}\sum_{i=1}^k |N(A_i)|+\sqrt{(\log d) \left(2^{d-1}-\sum_{i=1}^k |N(A_i)|\right)} \leq d(\log d)\left(\frac{2}{1+\gl}\right)^d.
$$
This completes the proof.
\qed

\medskip

Finally, we turn to (\ref{eq-scaling}). Note that the right-hand side of (\ref{eq-scaling}) is
$$
\Pr({\rm Poisson}(\grg_k)=c)
$$
where ${\rm Poisson}(\grg_k)$ is a Poisson random variable with parameter $\grg_k:=\frac{1}{2}e^{-k/2}$.

For each fixed $c \in {\mathbb N}$ we get a lower bound on the contribution to the partition function from those $I$ with $\min\{|I \cap \cE|, |I \cap \cO|\}=c$ by considering those which have $c$ $2$-components on $\cE$, each of size $1$, and have more than $\log d$ (say) vertices on $\cO$, and the same with $\cE$ and $\cO$ reversed. This gives a lower bound of
$$
\frac{2}{c!}\prod_{i=0}^{c-1} \left(2^{d-1}-id^2\right) \gl\left((1+\gl)^{2^{d-1}}-\sum_{i < \log d} \gl^i{2^{d-1}-cd \choose \log d}\right)
$$
which is at least
\begin{equation} \label{scaling-lb}
(2-o(1))(1+\gl)^{2^{d-1}}\frac{1}{c!}\left(\gld\right)^c.
\end{equation}
Recalling the discussion just before the proof of Lemma \ref{lem-small_no_of_comps_gl4}, we know that all but a vanishing part of $Z_\gl(Q_d)$ comes from $I$ with the smaller of $I \cap \cE$, $I \cap \cO$ consisting of no more than $\log d$ $2$-components of size $1$. So we get an upper bound on the contribution to the partition function from those $I$ with $\min\{|I \cap \cE|, |I \cap \cO|\}=c$ of
\begin{equation} \label{scaling-ub}
{2^{d-1} \choose c}\gl^c(1+\gl)^{2^{d-1}-cd} + o\left(Z_\gl(Q_d)\right) \leq (2+o(1))(1+\gl)^{2^{d-1}}\frac{1}{c!}\left(\gld\right)^c
\end{equation}
the inequality following from the fact that in this range of $\gl$,
$$
Z_\gl(Q_d) \sim 2(1+\gl)^{2^{d-1}}\exp\left\{\frac{\gl}{2}\left(\frac{2}{1+\gl}\right)\right\} = O\left((1+\gl)^{2^{d-1}}\right).
$$
Combining (\ref{scaling-lb}) and (\ref{scaling-ub}), and noting that for $\gl = 1+\frac{k+o(1)}{d}$, $\gld \sim \grg_k$, it follows that
$$
\Pr\left(\min\{|I \cap \cE|, |I \cap \cO|\}=c\right)  \sim  \Pr({\rm Poisson}(\grg_k)=c).
$$

\subsection{Proof of Theorem \ref{thm-structure}}
\label{subsec-proof_of_thm_structure}

The first statement follows immediately from the fact that (\ref{eq-thresh_gl1}) holds a.a.s. for $\gl$ satisfying (\ref{lambda1}). The second statement follows from our proof that (\ref{eq-scaling}) holds a.a.s. for $\gl$ satisfying (\ref{lambda1}), once we observe that the lower bound in (\ref{eq-scaling}) is obtained by only considering those independent sets for which ${\rm cl}(I_{\rm min}) \leq 1$. By a similar observation, our proof that (\ref{eq-thresh_gl3}) holds a.a.s. for $\gl$ satisfying (\ref{lambda3}) also proves the third statement.

For the fourth statement, note for $2^{1/m'}-1 \geq \gl > 2^{1/(m'+1)} - 1$ we satisfy (\ref{def-mgl}) with $m=m'$; the statement then follows immediately from Lemma \ref{lem-small_no_of_comps_gl4}.

\subsection{Proof of Theorem \ref{thm-influence}}
\label{subsec-proof_of_thm_influence}

Our approach is inspired by \cite{GalvinKahn}, in which Galvin and Kahn obtain a result of a similar flavour on the lattice ${\mathbb Z}^d$. We begin with \eqref{pt_odd}. Write
$$
\cJ = \{J \in \cI(Q_d):w \in J\} ~~~\mbox{and}~~~
\cI^\prime = \{I \in \cJ:u \in I\}.
$$
Further, write $\cI = \{I \in \cI^\prime:I \cap \cE~\mbox{small}\}$.
We need to bound
\begin{equation} \label{need0}
\frac{w_\gl(\cI^\prime)}{w_\gl(\cJ)} \leq
(1+\gl)^{-d(1-o(1))},
\end{equation}
where $w_\gl(*)=\sum_{I \in *} \gl^{|I|}$.
We will show
\begin{equation} \label{need}
\frac{w_\gl(\cI)}{w_\gl(\cJ)} \leq
(1+\gl)^{-d(1-o(1))}.
\end{equation}
The same argument will show
$$
\frac{w_\gl(\{I \in \cI^\prime:I \cap \cO~\mbox{small}\})}{w_\gl(\{J \in \cI(Q_d):u \in J\})} \leq (1+\gl)^{-d(1-o(1))}.
$$
Combining this with (\ref{need}) we get
(\ref{need0}), noting that for any $I \in
\cI(Q_d)$, either $I \cap \cE$ or $I \cap \cO$ small, and that by symmetry $w_\gl(\{J \in \cI(Q_d):u \in J\})=w_\gl(\cJ)$.

We will obtain (\ref{need}) by producing, for each $I \in \cI$, a
set $\varphi(I) \subseteq \cJ$, as well as a map $\nu:\cI \times \cJ
\rightarrow {\mathbb R}$ supported on pairs $(I,J)$ with $J \in
\varphi(I)$ and satisfying
\begin{equation} \label{nuout}
\sum_{J \in \varphi(I)} \nu(I,J) = 1
\end{equation}
for each $I \in \cI$ and
\begin{equation} \label{nuin}
\sum_{I \in \varphi^{-1}(J)} \gl^{|I|-|J|}\nu(I,J) \leq
(1+\gl)^{-d(1-o(1))}
\end{equation}
for each $J \in \cJ$. It is not difficult to see that the existence of such a $\varphi$
and $\nu$ satisfying (\ref{nuout}) and (\ref{nuin}) gives
(\ref{need}).

We produce $\varphi$ as follows. Given $I \in \cI$, write $W(I)$ for
the $2$-component of $I \cap \cE$ containing $u$. Set
$$
\cW(a,g) = \{W \subseteq \cE: |W|=a,~|N(W)|=g,~u \in
W,~W~\mbox{small and $2$-linked}\}
$$
and
$$
\cI(a,g) = \{I \in \cI : W(I) \in \cW(a,g)\}.
$$
Set $I^\prime = I \setminus W$. Note that $N(W(I)) \cap I = \emptyset$ and
$N(W(I))$ is not adjacent to anything in $I'$ (if it was, then $W(I)$
would not be the $2$-component of $u$ in $I \cap \cE$). We may
therefore add any subset of $N(W(I))$ to $I'$ and still have an
independent set. Set
$$
\varphi(I) = \{I^\prime \cup S: S \subseteq N(W(I))\}.
$$
We have just observed that indeed $\varphi(I) \subseteq \cJ$.

For each $J \in \varphi(I)$ write $S(J)$ for $J \setminus I'$ and
set
$$
\nu(I,J) = \frac{\gl^{S(J)}}{(1+\gl)^{|N(W(I))|}}
~\left(=\frac{\gl^{|J|-|I|+|W(I)|}}{(1+\gl)^{|N(W(I))|}}\right).
$$
Since $S(J)$ runs over all subsets of $N(W(I))$ it is clear that
(\ref{nuout}) holds. To see that (\ref{nuin}) holds, observe that
for fixed $J \in \cJ$ we have
\begin{eqnarray}
\sum_{I \in \varphi^{-1}(J)} \gl^{|I|-|J|}\nu(I,J) & = & \sum_{I
\in \varphi^{-1}(J)} \gl^{|W(I)|}(1+\gl)^{-N(W(I))} \nonumber \\
& \leq &  \sum_{a,g,W} ~~~\sum_{I \in \varphi^{-1}(J), ~I \in
\cI(a,g),~W(I)=W} \gl^{a}(1+\gl)^{-g} \nonumber \\
& \leq & \sum_{a,g} \sum_{W \in \cW(a,g)} \gl^{a}(1+\gl)^{-g}
\label{unique_reconstruction}.
\end{eqnarray}
The main point here is (\ref{unique_reconstruction}), which follows
from the fact that for each $W \in \cW(a,g)$ and $J \in \cJ$ there
is at most one $I \in \cI$ such that $I \in \cI(a,g)$, $W(I)=W$ and
$I \in \varphi^{-1}(J)$.

For each $g>d^4$ we have $\cW(a,g) \subseteq \cup_{a' \leq g}
\cG(a',g)$, and so, using Lemma
\ref{lem-mainsaplemma_weighted} for (\ref{using_saplemma}),
\begin{eqnarray}
\sum_{a,~g>d^4,~W \in \cW(a,g)} \gl^{a}(1+\gl)^{-g} & \leq &
\sum_{g>d^4,~a' \leq g} (1+\gl)^{-g}\sum_{A \in
\cG(a',g,u)}\gl^{|A|}
\nonumber \\
& \leq & 2^d \sum_{g>d^4,~a' \leq g}
(1+\gl)^{-\frac{c'(g-a')\log d}{d^{2/3}}}. \label{using_saplemma}
\end{eqnarray}
By Lemma \ref{lem-main-cube-iso} we have $g-a'=\Omega(d^{7/2})$ in the range
$g>d^4$ and we have at most $2^d$ choices for each of $a'$ and $g$
and so
\begin{equation} \label{case_large}
\sum_{a,g>d^4,W \in \cW(a,g)} \gl^{a}(1+\gl)^{-g} \leq  (1+\gl)^{-d(1-o(1))}.
\end{equation}

For $g \leq d^4$ we have $|\cW(a,g)| \leq 2^{O(a \log d)} \leq
2^{O(g\log d/d)}$ and so
\begin{eqnarray}
\sum_{a,~g\leq d^4,~W \in \cW(a,g)} \gl^{a}(1+\gl)^{-g} & \leq &
\sum_{a,~g\leq d^4}
(1+\gl)^{-g} 2^{O(g\log d/d)} (1+\gl)^{O(g/\log d)} \nonumber \\
& \leq & \sum_{g\geq d} (1+\gl)^{\left\{O\left(\frac{\log d}{\log
(1+\gl)}\right)-g+O\left(\frac{g \log d}{d \log
(1+\gl)}\right)+\frac{g}{\log d}\right\}}. \nonumber
\\
& \leq & (1+\gl)^{-d(1-o(1))} \label{case_small}.
\end{eqnarray}
Combining \eqref{case_large} with \eqref{case_small} we obtain
\eqref{need} and so \eqref{need0} and \eqref{pt_odd}.

We obtain \eqref{pt_even} from \eqref{pt_odd} easily. Conditioned on
$v \in I$, the probability that a particular neighbour
of $u$ is in $I$ is, by \eqref{pt_odd}, at most
$(1+\gl)^{-d(1-o(1))}$, and so the probability that none of the $d$
neighbours of $u$ are in $I$ is at least
$1-d(1+\gl)^{-d(1-o(1))} = 1-o(1)$. The probability that $u$ is in
$I$ is at least the probability that it is in $I$ conditioned on none of neighbours being in $I$ times
the probability that none of neighbours are in $I$, and so
is at least $(1-o(1))\gl/(1+\gl)$.

\section{Proof of Lemma \ref{lem-mainsaplemma_weighted}} \label{sec-sap_proof}

There are three steps to the proof. In the first step (Lemma \ref{lem-firstapprox}) we associated to each $A \in \cG(a,g)$ a pair $(F^\star, S^\star)$ that approximates $A$ in the sense that $F^\star \subseteq N(A)$, $S^\star \supseteq [A]$ and both of $|N(A)\setminus F^\star|, |[A]\setminus S^\star| \leq x$ hold for some suitably small $x$, and we bound the size of $\cA_1$, the set of all pairs $(F^\star, S^\star)$ that arise as we run over $A \in \cG(a,g)$ (the bound, of course, depending on $x$ as well as $a$ and $g$). This first step is the most involved of the three, and our presentation of it is based closely on Sapozhenko's original treatment \cite{Sapozhenko}.

The first step may be thought of as a partitioning of $\cG(a,g)$, with the $|\cA_1|$ many partition classes indexed by pairs $(F^\star, S^\star)$. The second step (Lemma \ref{lem-secondapprox}) focuses on the individual partition classes: to each $(F^\star, S^\star)$ and $A$ in the class indexed by $(F^\star, S^\star)$ we associated a pair $(F,S)$ that approximates $A$ in the sense that $F \subseteq N(A)$, $S \supseteq [A]$ and $|S| \leq |F| + y$ for some suitably small $y$, and we bound (uniformly in $(F^\star, S^\star)$) the size of $\cA_2$, the set of all pairs $(F, S)$ that arise as we run over $A$ in the class indexed by $(F^\star, S^\star)$ (the bound depending on $y$). This second step essentially appears in work of Galvin and Kahn \cite{GalvinKahn} (with a proof also adapted from \cite{Sapozhenko}), and here we only show how the conclusion of \cite[Lemma 2.17]{GalvinKahn} almost immediately yields our desired conclusion.

In the third step (Lemma \ref{lem-reconstruction}) we focus on a particular $(F,S)$ and bound (uniformly in $(F,S)$) the sum of the $\gl^{|A|}$'s over all those $A \in \cG(a,g)$ for which it holds that $F \subseteq N(A)$, $S \supseteq [A]$ and $|S| \leq |F| + y$. This comes directly from work of Galvin and Tetali \cite{GalvinTetali2} and so we do not give the proof here.

The steps together give
$|\cG(a,g)| \leq B_\gl|\cA_1||\cA_2|$
where $B_\gl$ is the bound from the third step. We present the steps in more generality than we need, since this adds nothing to the complexity of the proofs.

\begin{lemma} \label{lem-firstapprox}
Let $\gS$ be a $d$-regular bipartite graph with bipartition classes $X$ and $Y$.
Let $\cG = \{A \subseteq X~\mbox{$2$-linked}: |[A]|=a, |N(A)|=g\}$.
Fix $1 \leq \varphi \leq d-1$. Let
$$
m_\varphi = \min \left\{N(K):y \in Y, K \subseteq N(y), |K|>\varphi \right\}.
$$
Let $C>0$ be any constant such that $C\log d /(\varphi d) < 1$. Set $t=g-a$. There is a family ${\cal A}_1 \subseteq 2^Y \times 2^X$ with
\begin{eqnarray*}
|{\cal A}_1| & \leq & |Y|\exp\left\{\frac{78gC \log^2 d}{\varphi d} + \frac{78g \log d}{d^{Cm_\varphi/(\varphi d)}} + \frac{78t\log^2 d}{d-\varphi}\right\} \times \\
& & ~~~~~~~~~~~~~~~~~~~~~~~~~~~~{\frac{3gC \log d}{\varphi} \choose \leq \frac{3tC \log d}{\varphi}}{dg \choose \leq dt/(\varphi(d-\varphi))}
\end{eqnarray*}
(where ${n \choose \leq k}$ is shorthand for $\sum_{i \leq k} {n \choose i}$) and a map $\pi_1:\cG \rightarrow \cA_1$ for which $\pi_1(A):=(F^\star, S^\star)$ satisfies $F^\star \subseteq N(A)$, $S^\star \supseteq [A]$, $|N(A) \setminus F^\star| \leq td/(d-\varphi)$ and $|S^\star \setminus [A]| \leq  td/(d-\varphi)$.
\end{lemma}

\begin{lemma} \label{lem-secondapprox}
Let $\gS$ and $\cG$ be as in Lemma \ref{lem-firstapprox}. Let $(F^\star,S^\star) \in 2^Y \times 2^X$ and $x > 0$ be given. Let
$$
\cG' = \{A \in \cG: F^\star \subseteq N(A), S^\star \supseteq [A], |N(A) \setminus F^\star| \leq x~\mbox{and}~|S^\star \setminus [A]| \leq x\}.
$$
There is a constant $c>0$, a family ${\cal A}_2 \subseteq 2^Y \times 2^X$ with
$$
|{\cal A}_2| \leq \exp \left\{\frac{cx}{d}+\frac{ct\log d}{\psi}\right\}
$$
and a map $\pi_2:\cG' \rightarrow \cA_2$ for which $\pi_2(A):=(F, S)$ satisfies $F \subseteq N(A)$, $S \supseteq [A]$ and
\begin{equation} \label{boundingsbyf}
|S| \leq |F| + 2t\psi/(d-\psi).
\end{equation}
\end{lemma}

\begin{lemma} \label{lem-reconstruction} [See \cite[Lemma 3.4]{GalvinTetali2}]
Let $\gS$ and $\cG$ be as in Lemma \ref{lem-firstapprox}.
Let $\psi$ and $\grg$ satisfy $1 \leq \psi \leq d/2$ and $1 \geq \grg >
\frac{-2\psi}{d-\psi}$. Fix $(F,S) \in 2^Y \times 2^X$
satisfying (\ref{boundingsbyf}). We have
$$
\sum \gl^{|A|} \leq \max \left\{(1+\gl)^{g-\grg t}, {3dg \choose
\leq \frac{2t\psi}{d-\psi}+\grg t}(1+\gl)^{g-t}\right\}
$$
where the sum is over all $A \in \cG$ satisfying $F
\subseteq N(A)$ and $S \supseteq [A]$.
\end{lemma}

Before turning to the proofs, we put them together in the case $\gS=Q_d$. We set $\varphi=d/2$ (which choice allows us to take $x=2t$ in Lemma \ref{lem-secondapprox}) and $\psi=d^{2/3}$. By Lemma \ref{lem_easy-cube-iso} we have $m_\varphi \geq d^2/(2C_{\rm iso}$), and (for large enough $d$) we may set $C=2C_{\rm iso}$. Using $t \geq \Omega(g/\sqrt{d})$ (from Lemma \ref{lem-main-cube-iso}) and the basic binomial estimate
\begin{equation} \label{basic-bin}
{n \choose \leq k}  \leq  \exp\left\{(1+o(1))\left(k\log\frac{n}{k}\right)\right\}
\end{equation}
for $k=o(n)$, the first two lemmas combine to give
$$
|\cA_1||\cA_2| \leq 2^d \exp \left\{O\left(\frac{t\log d}{d^{2/3}}\right)\right\}.
$$
We now take
$$
\grg = \frac{\log (1+\gl) - \frac{6\psi \log d}{d-\psi}}{\log(1+\gl)+3\log d} \geq \frac{\frac{c}{3}-3}{d^{1/3}},
$$
with the inequality valid for $\gl>\bestgl$.
By our choices of $\psi$ and $\grg$ we have
$$
\frac{g}{d} \leq \frac{2t\psi}{d-\psi}+\grg t \leq 3g
$$
and so
\begin{eqnarray*}
{3dg \choose \leq \frac{2t\psi}{d-\psi}+\grg t} & \leq & \exp\left\{(1+o(1))\left(\frac{2t\psi}{d-\psi}+\grg t\right)\log \left(\frac{3dg}{\frac{2t\psi}{d-\psi}+\grg t}\right)\right\} \\
& \leq & \exp\left\{3\left(\frac{2t\psi}{d-\psi}+\grg t\right)\log d\right\} \\
& = & (1+\gl)^{(1-\grg)t}
\end{eqnarray*}
with the first inequality using (\ref{basic-bin}) and the equality following from the definition of $\grg$. It follows that
\begin{eqnarray*}
\max \left\{(1+\gl)^{g-\grg t}, {3dg \choose
\leq \frac{2t\psi}{d-\psi}+\grg t}(1+\gl)^{g-t}\right\} & \leq & (1+\gl)^{g-\grg t}
\end{eqnarray*}
so that
\begin{eqnarray*}
|\cG(a,g)| & \leq & 2^d(1+\gl)^g \exp\left\{-\grg t \log (1+\gl) + O\left(\frac{t\log d}{d^{2/3}}\right)\right\} \\
& \leq & 2^d(1+\gl)^g \exp\left\{-\frac{c't\log d}{d^{2/3}}\right\}
\end{eqnarray*}
for some $c'>0$ (as long as $c>0$ is suitably large), as claimed.

\medskip

\noindent {\em Proof of Lemma \ref{lem-firstapprox}}: Fix $A \in \cG$ and set
$$
N(A)^{\varphi} = \{y \in N(A):d_{[A]}(y)>\varphi\}
$$
(where for any $K \subseteq V$, $d_{K}(y) := |N(y)\cap K|$).
We begin by describing the construction of an $F'$ which satisfies $N(A)^\varphi \subseteq F' \subseteq N(A)$ and $N(F') \supseteq [A]$. Since each vertex in $N(A) \setminus F'$ is in $N(A) \setminus N(A)^{\varphi}$ and so contributes at
least $d-\varphi$ edges to $\nabla(N(A), X\setminus [A])$, a set of size
$gd-ad=td$, such a set satisfies $|N(A) \setminus F'| \leq td/(d-\varphi)$. (Here and throughout we use $\nabla(A,B)$ to indicate the set of edges with one endpoint in $A$ and the other in $B$.)

Set $p=C \log d /(\varphi d)$.
Construct a random subset $\tilde{T}$ of $N(A)$ by
putting each $y \in N(A)$ in $\tilde{T}$ with probability $p$, these choices
made independently. We have
\begin{equation} \label{towardst0}
{\bf E}(|\tilde{T}|)=gp
\end{equation}
and since $|\nabla(N(A),X\setminus [A])| = td$,
\begin{equation} \label{towardsomega}
{\bf E}(|\nabla(\tilde{T},X\setminus [A])|) = tdp.
\end{equation}
For $y \in N(A)^{\varphi}$ we have $|N(N_{[A]}(\{y\}))| \geq m_\varphi$
and so
\begin{eqnarray}
{\bf E}(|N(A)^{\varphi} \setminus N(N_{[A]}(\tilde{T}))|) & = & \sum_{y \in
N(A)^{\varphi}} \Pr(y \not \in
N(N_{[A]}(\tilde{T}))) \nonumber \\
& = & \sum_{y \in N(A)^{\varphi}} \Pr(N(N_{[A]}(\{y\})) \cap \tilde{T} =
\emptyset) \nonumber \\
& \leq & g(1-p)^{m_\varphi} \nonumber \\
& < & g\exp\left\{-pm_\varphi\right\}. \label{towardst0prime}
\end{eqnarray}
Combining (\ref{towardst0}), (\ref{towardsomega}) and
(\ref{towardst0prime}) and using Markov's inequality we find that
there is at least one $T_0 \subseteq N(A)$ satisfying
\begin{equation} \label{t0}
|T_0|\leq \frac{3C g \log d}{\varphi d}
\end{equation}
\begin{equation} \label{omega}
|\Omega| \leq \frac{3C td \log d}{\varphi d}
\end{equation}
where $\Omega := \nabla(T_0,X\setminus [A])$ and
\begin{equation} \label{t0prime}
|N(A)^{\varphi} \setminus N(N_{[A]}(T_0))| \leq \frac{3g}{d^{Cm_\varphi/(\varphi d)}}.
\end{equation}
Choose one such $T_0$ and set $T_0':=N(A)^{\varphi} \setminus N(N_{[A]}(T_0))$. Setting $L=N(N_{[A]}(T_0)) \cup T_0'$, we have $L \supseteq N(A)^{\varphi}$. Let $T_1 \subseteq N(A) \setminus L$ be a cover of minimum size of $[A]
\setminus N(L)$ in the graph induced by $(N(A) \setminus L) \cup ([A]
\setminus N(L))$. Set $F'=L \cup T_1$. By construction, $F'$ satisfies $N(A)^\varphi \subseteq F' \subseteq N(A)$ and $N(F') \supseteq [A]$.

Before estimating how many sets $F'$ might be produced in this way
as we run over $A\in {\cal G}$, we make some observations about
the sets described above.

First, note that by Lemma \ref{nearbysets}, $F'$ is $4$-linked ($[A]$
is $2$-linked, every $x\in [A]$ is at distance $1$ from $F'$ and every
$y \in F'$ is at distance $1$ from $[A]$) and so $T=T_0 \cup T_0' \cup T_1$ is $8$-linked (every $y
\in T$ is at distance $2$ from $F'$ and every $y \in
F'$ is at distance $2$ from $T$).

Next, note that $F'$ is completely
determined by the tuple $(T_0, T_0', T_1, \Omega)$, since $T_0$ and $\Omega$ together determine $N(N_{[A]}(T_0))$.

The sizes of $T_0$, $T_0'$ and $\Omega$ are bounded by
(\ref{t0}), (\ref{t0prime}) and (\ref{omega}), respectively.
To bound $|T_1|$, note that as previously observed $|N(A) \setminus L| \leq td/(d-\varphi)$,
$d_{[A] \setminus N(L)}(u) \leq d$ for
each $u \in G \setminus L$, and $d_{G \setminus L}(v) = d$ for each
$v \in [A] \setminus N(L)$. So by Lemma \ref{lovaszstein},
$|T_1|
\leq (t/(d-\varphi)) (1+\log d) \leq 3t\log d/(d-\varphi)$.

Combining these observations, we get that $T$ is an $8$-linked
subset of $Y$ with
$$
|T|\leq \frac{3gC \log d}{\varphi d} + \frac{3g}{d^{Cm_\varphi/(\varphi d)}} + \frac{3t\log d}{d-\varphi} =: T_{\rm bound}.
$$
By Corollary \ref{Tree} there are
$|Y|\exp\left\{24T_{\rm bound}\log d\right\}$
possible choices for $T$. Once $T$ has been chosen, there are at most $2^{T_{\rm bound}}$
choices for $T_0 \subseteq T$, at most $2^{T_{\rm bound}}$ choices for $T_1 \subseteq T$ and at most
$$
{\frac{3gC \log d}{\varphi} \choose \leq \frac{3tC \log d}{\varphi}}
$$
choices for $\Omega$. So the number of choices for $(T_0, T_0', T_1, \Omega)$ is at most
\begin{equation} \label{int17}
|Y| \exp\left\{26T_{\rm bound}\log d\right\} {\frac{3gC \log d}{\varphi} \choose \leq \frac{3tC \log d}{\varphi}}
\end{equation}

We now describe an algorithmic procedure which produces $F^\star$ from $F'$, and also produces $S^\star$ (again, for a fixed $A$).
If $\{u \in [A]: d_{N(A) \setminus F'}(u)
>\varphi \}
\neq \emptyset$, pick the smallest (with respect to some fixed ordering of the vertices of $\gS$) $u$ in
this set and update $F'$ by $F' \longleftarrow F' \cup N(u)$.
Repeat this until $\{u \in [A]: d_{N(A) \setminus F'}(u) > \varphi\} =
\emptyset$. Then set $F^\star=F'$ and $S^\star=\{u \in X:d_{F^\star}(u)
\geq d-\varphi\}$.

Observe that $F^\star$ thus constructed inherits the properties $F^\star \subseteq N(A)$ and $|N(A)\setminus F^\star| \leq td/(d-\varphi)$ from $F'$, since we obtain $F^\star$ from $F'$ by adding vertices of $N(A)$. We also have $S^\star \supseteq [A]$, since otherwise the algorithm would not have terminated. Since each vertex of $S^\star \setminus [A]$ contributes at least $d-\varphi$ vertices to $\nabla(N(A),X \setminus [A])$, a set of size $td$, we have $|S^\star \setminus [A]|\leq td/(d-\varphi)$. Finally, the algorithm is determined by the selection of at most $dt/(\varphi(d-\varphi))$ vertices (each iteration removes at least $\varphi$ vertices from $N(A)\setminus F'$, a set of initial size at most $td/(d-\varphi)$). These vertices come from $[A]$ which is contained in $N(F')$, a set of size at most $dg$. So the total number of possibilities for $(F^\star, S^\star)$ for each $F'$ is at most
$$
{dg \choose \leq dt/(\varphi(d-\varphi))}.
$$
Combining this with (\ref{int17}), we obtain the claimed bound on $|\cA_1|$.
\qed

\medskip

\noindent {\em Proof of Lemma \ref{lem-secondapprox}}: An almost identical statement appears in \cite[Lemma 2.17]{GalvinKahn}, the difference being that (\ref{boundingsbyf}) is replaced by the two conditions $d_F(u) \geq d-\psi$ for all $u \in S$ and $d_{X \setminus S}(v) \geq
d-\psi$ for all $v \in Y \setminus F$. (The proof essentially repeats the algorithmic procedure described at the end of the proof of Lemma \ref{lem-firstapprox}, with $\varphi$ replaced by $\psi$). But these two degree conditions imply (\ref{boundingsbyf}).
Indeed, observe that $|\nabla(S,G)|$ is bounded above
by $d|F| + \psi|N(A) \setminus F|$ and below by $d|[A]| + (d-\psi)|S
\setminus [A]| = d|S| - \psi|S \setminus [A]|$, giving
$$
|S| \leq |F| + \psi|(N(A) \setminus F) \cup (S \setminus [A])|/d,
$$
and that each $u \in (N(A) \setminus F) \cup (S \setminus [A])$
contributes at least $d-\psi$ edges to $\nabla(N(A),X\setminus [A])$, a
set of size $td$, giving
$$
|(N(A) \setminus F) \cup (S \setminus A)| \leq 2td/(d-\psi).
$$
These two observations together give (\ref{boundingsbyf}). \qed

\bigskip

{\em Acknowledgements}: Part of this work was carried while the author was a participant in the programme on Combinatorics and Statistical Mechanics at the Isaac Newton Institute for Mathematical Sciences, University of Cambridge, in spring 2008. The author thanks the Institute and the programme organizers for the support provided. The author also thanks the referee for a careful reading of the manuscript and for many helpful comments.

\end{document}